\documentclass[smallextended]{svjour3}       
\smartqed  
\usepackage{graphicx}
\usepackage{amsmath}
\usepackage{amssymb, comment}
\usepackage{mathrsfs}
\usepackage{multirow}
\numberwithin{equation}{section}

\def \E {\mathbb E}

\def \P {\mathbb{P}}

\renewcommand{\proofname}{\noindent {\bf Proof}}

\title{Heavy-traffic analysis of the $M^X/\text{semi-Markov}/1$ queue}

\begin{document}

\author{Abhishek \and
        Marko~Boon \and
        Rudesindo~N\'u\~nez-Queija
        }

\institute{Abhishek $\&$  Rudesindo~N\'u\~nez~Queija \at
           Korteweg-de Vries Institute for Mathematics, University of Amsterdam, Amsterdam, The Netherlands\\
           \email{\{Abhishek, nunezqueija\}@uva.nl}
         \and
        Marko~Boon \at
        Department of Mathematics and Computer Science, Eindhoven University of Technology, Eindhoven, The Netherlands \\
        \email{m.a.a.boon@tue.nl}
       }

\date{\today}
\maketitle
\begin{abstract}
In this paper we analyze a single server queue with batch arrivals and semi-Markovian service times. We also include the feature that the first service of each busy period might have a different distribution than subsequent service times. Our generating function based approach allows us to determine the heavy traffic limit of the scaled queue-length distribution.
It turns out that this distribution converges to an exponential distribution. Nonsurprisingly, the exceptional first service does not influence this limiting distribution. We identify a sufficient and necessary condition under which the dependence between successive service times disappears in the limit, which we illustrate in a numerical example.

\end{abstract}

\section{Introduction}
In many systems, successive service times of customers are not independent. The service type of a customer may depend on the type and the service duration of the preceding customer. Queueing systems with correlated service durations arise in many applications: logistics, production/inventory systems, computer and telecommunication networks. The model considered in this paper, is specifically  motivated by its application to a road traffic setting, in which a stream of vehicles on a minor road merges with, or crosses a stream on a main road at an unsignalized intersection. The queueing model of vehicles on the minor road is known to be a single server queue (possibly with batch arrivals) with semi-Markovian service times and a different service time distribution for vehicles that arrive when no queue is present (see \cite{AbhishekMMOR,AbhishekWaitingTimes}). The dependence between successive service times (i.e. the time to wait for a sufficiently large gap and crossing the road) is caused by either platoon forming on the major road, or by the fact that the remaining part of a gap between successive vehicles on the major road may be used by the following vehicle.

In this paper, we consider a single-server queue with batch arrivals and correlated service times. The correlations are modeled with different service types, which form a Markov chain that itself depends on the sequence of service lengths. In addition, the first customer in a busy period has a different service time distribution  than regular customers served in the busy period, which was firstly introduced in the framework of the $M/G/1$ queueing model by Welch \cite{welch} and by Yeo \cite{yeo}.

Queues with  correlated service times have been studied for many years \cite{QUESTA2017,cin_s,gaver,neuts66,neuts77a,neuts77b}. One of the first studies for Markov-modulated single-server queueing systems in heavy traffic (HT) was by Burman and Smith \cite{burman}, who study the mean delay and the mean number in queue in a single-server system in both light-traffic and heavy-traffic regimes, where customers arrive according to a nonhomogeneous Poisson process with rate equal to a function of the state of an independent Markov process. In their model, service times are independent and identically distributed. Later, G. Falin and A. Falin \cite{falin} suggest another approach to analyze the same queueing model, which is based on certain  ‘semi-explicit’ formulas for the stationary distribution of the virtual waiting time and its mean value under heavy traffic. Dimitrov \cite{dimitrov} applies the same approach to a single-server queueing system with arrival rate and service time depending on the state of Markov chain at an arrival epoch, and shows that the distribution of the scaled stationary  virtual waiting time is exponential under a HT scaling. Several other authors \cite{asmussen,thorsdottir} also study Markov-modulated $M/G/1$-type queueing systems in heavy traffic. However, we are not aware of any prior  work analyzing the $M^X/\text{SM}/1$ queue with exceptional first service under a HT scaling.

The current model is a slight extension of that in \cite{QUESTA2017}. We allow the service duration of a customer arriving into an empty system to have a distribution that differs from the service-time distributions of other customers. For the stationary analysis of the model this requires minor adaptations of that in \cite{QUESTA2017}. In addition, we investigate the stationary distribution in the heavy-traffic regime.

The remainder of this paper is organized as follows. In Section \ref{description}, we present the description of the queueing model. In Section \ref{stationary queue length}, we first determine the stationary probability generating function of the queue length of the system  at the departure time of a customer. Subsequently, we use that result to derive the generating functions of the stationary queue length at an arbitrary time, at batch arrival instants, and at customer arrival instants. Using these results, we obtain  the heavy-traffic distribution of the scaled stationary queue length in Section \ref{heavy-traffc}. In Section \ref{numerical_results}, a numerical example is presented to demonstrate the impact of the correlated service times on the queue length distribution in the heavy-traffic regime.

\section{Model description}\label{description}
We consider a single-server queuing system. Customers arrive in batches  at the system according to a Poisson process  with rate $\lambda$. The arriving batch size is denoted by the random variable $B$, with generating function $B(z)$, for $|z|\leq 1$ (zero-sized batches are not allowed, i.e. $B\geq 1$). Customers are served individually, and the first customer in a busy period has a different service time distribution  than regular customers served in the busy period.  There are $N$ types of customers, which we number $1,2,\dots,N$. Denote by $J_n$ the type of the $n$th customer and $G^{(n)}$ its service time, $n=1,2,\dots$. The type of a customer is only determined at the moment its service begins. More specifically, the type of the $n$th customer depends on the type, and on  the service duration of the $(n-1)$th customer, as well as on whether the queue is empty at the departure time of the $(n-1)$th customer. We introduce, for $i=1,2,\dots,N$,
\begin{align}
\tilde{G}_{ij}(s)&=\E[e^{-sG^{(n)}}1_{\{J_{n+1}=j\}}|J_n=i, X_{n-1}\geq 1],\label{G_{ij}(s)}\\
\tilde{G}^*_{ij}(s)&=\E[e^{-sG^{(n)}}1_{\{J_{n+1}=j\}}|J_n=i, X_{n-1}=0]\label{G*_{ij}(s)},
\end{align} where $X_{n-1}$ is the number of customers in the system immediately after the departure of the $(n-1)$th customer. \\

In particular, for $i,j=1,2,\dots,N$, we define
\begin{align}
P_{ij}=\tilde{G}_{ij}(0)=\P(J_{n+1}=j|J_n=i,X_{n-1}\geq 1),\label{p_ij}\\
P^*_{ij}=\tilde{G}^*_{ij}(0)=\P(J_{n+1}=j|J_n=i,X_{n-1}= 0).\label{p^*_ij}
\end{align}
In the literature, the service process considered in this paper is referred to as a semi-Markov (SM) process, and thus the queuing system is referred to as the $M^X/SM/1$. In fact, in the gap acceptance literature, the single server queue with exceptional first service is commonly referred to as the $M/G2/1$ queue, a term seemingly introduced by Daganzo \cite{daganzo1977}, which motivates us to denote this model (with batch arrivals and exceptional first service) as the $M^X/SM2/1$ queue.\\

We assume that $P=[P_{ij}]_{i,j\in\{1,2,\dots,N\}}$ is the transition probability matrix of an irreducible discrete time Markov chain, with stationary distribution $\pi=(\pi_1,\pi_2,\dots,\pi_N)$ such that
\begin{align}\label{rel: piandP}
  \pi P=\pi.
\end{align}
For intuition we may think of $\pi$ as the conditional equilibrium distribution of $J_n$ in case the queue would never empty.
Using Cramer's rule with the normalizing equation $\sum_{i=1}^{N}\pi_i=1$, the solutions of the system of equations \eqref{rel: piandP} are given by
\begin{align}\label{pi_i}
 \pi_i=\frac{d_i}{d},
\end{align}
where $d=\sum_{i=1}^Nd_i,$ and $d_i$ is the cofactor of the entry in the $i$-th row and the first column of the matrix  $(I-P)$, which is given by
\begin{align}
 d_1=&\begin{vmatrix}
1-P_{22} & -P_{23} & \dots & -P_{2N}\\
-P_{32} & 1-P_{33} & \dots & -P_{3N}\\
\vdots & \vdots & \ddots & \vdots \\
-P_{N2} & -P_{N3} & \dots & 1-P_{NN}
\end{vmatrix}, \label{d_1}\\
 d_i=&(-1)^{i+1}\begin{vmatrix}
-P_{12} & -P_{13} & \dots & -P_{1N}\\
\vdots & \vdots & \ddots & \vdots \\
-P_{i-12} & -P_{i-13} & \dots & -P_{i-1N}\\
-P_{i+12} & -P_{i+13} & \dots & -P_{i+1N}\\
\vdots & \vdots & \ddots & \vdots \\
-P_{N2} & -P_{N3} & \dots & 1-P_{NN}
\end{vmatrix},\quad i=2,3,\dots,N-1, \label{d_i}\\
 d_N=&(-1)^{N+1}\begin{vmatrix}
-P_{12} & -P_{13} & \dots & -P_{1N}\\
1-P_{22} & -P_{23} & \dots & -P_{2N}\\
\vdots & \vdots & \ddots & \vdots \\
-P_{N-12} & -P_{N-13} & \dots & -P_{N-1N}
\end{vmatrix}. \label{d_N}
\end{align}
In the next section, to study the queue length distribution at departure times of customers, we denote by $A_{n}$ the number of arrivals during the service time of the $n$th customer (counting the individual customers inside the batches). We introduce, for $i=1,2,\dots,N$, \begin{align}
  A_i(z)&=\sum_{j=1}^{N}A_{ij}(z),\label{A_i(z)}\\
 A^{*}_i(z)&=\sum_{j=1}^{N}A^{*}_{ij}(z),\label{A^*_i(z)}
\end{align}  with
\begin{align}
A_{ij}(z)=\E[z^{A_{n}}1_{\{J_{n+1}=j\}}|J_n=i,X_{n-1}\geq 1], \label{A_ij(z)}\\
A^{*}_{ij}(z)=\E[z^{A_{n}}1_{\{J_{n+1}=j\}}|J_n=i,X_{n-1}=0]. \label{A_ijstar(z)}
\end{align}
Let us define
\begin{align}\label{rho}
\rho&=\sum_{i=1}^{N}\pi_i\alpha_i,
\end{align}
where
\begin{align}
\alpha_i&=\sum_{j=1}^N\alpha_{ij}, \label{alph_i}
\end{align} with
\begin{align}\label{alpha_ij}
\alpha_{ij}=\E[A_n 1_{\{J_{n+1}=j\}}|J_n=i,X_{n-1}\geq 1].
\end{align}
Intuitively, we can think of $\rho$ as being the expected number of arrivals during a service time if the process $(J_n,X_{n-1})$ would never hit the level $X_{n-1}=0$. Introducing some further notations:
\begin{align}
\alpha^{*}_i=\sum_{j=1}^N\alpha^{*}_{ij}, \label{alph_i*}
\end{align} with
\begin{align}\label{alpha_ij*}
\alpha^*_{ij}=\E[A^*_n 1_{\{J_{n+1}=j\}}|J_n=i,X_{n-1}=0].
\end{align}


Note that the number of arrivals during the service time of a customer is a batch Poisson process. Therefore, we can write the following relations:
\begin{align}
  A_{ij}(z)=&\tilde{G}_{ij}(\lambda(1-B(z))),\label{rel_A_G} \\
  A^{*}_{ij}(z)=&\tilde{G}^{*}_{ij}(\lambda(1-B(z))), \quad \text{ for } i,j=1,2,\dots,N. \label{rel_A*_G*}
\end{align}

To derive the stability condition for our model we use the results from \cite{QUESTA2017}. Note that the dynamics in the current model only differs from that in \cite{QUESTA2017} when the queue length is zero. More specifically, the two processes have identical transition rates, except in a finite number of states. This implies that the two processes are either both positive recurrent, both null recurrent or both transient. The condition for stability reads $\rho<1$, in accordance with \cite{QUESTA2017}, and similarly, both processes are null recurrent if $\rho=1$. Hence, if we modify the parameters such that $\rho\uparrow1$, the processes move from positive recurrence to null recurrence. In particular, $\P[X=0]>0$ if $\rho<1$ and $\P[X=0]\rightarrow0$ as $\rho\uparrow1$.


\section{ Stationary queue length analysis}\label{stationary queue length}
In this section, we shall first determine the steady-state joint distribution of the number of customers in the system immediately after a departure, and the type of the next customer to be served. Subsequently, we will use this result to derive the generating functions of the stationary number of customers at an arbitrary time, at batch arrival instants, and at customer arrival instants.

\subsection{ Stationary queue length analysis: departure epochs}
Starting-point of the analysis is the following recurrence relation:

\begin{align}
X_n = \left\{
\begin{array}{l l}
\ X_{n-1}-1+A_n & \quad \text{if $X_{n-1} \geq 1$ }\\
A_n +B_n-1& \quad \text{if $X_{n-1} =0$}
\end{array} \right., ~~~ n=1,2,3,\dots,
\label{recurA}
\end{align} where $X_n$ is the number of customers at the departure times of the $n$th customer and $B_n$ is the size of the batch in which $n$th customer arrived, with generating function $B(z)$, for $|z|\leq 1$. Due to dependent successive service times, $X_n$ here is not a Markov chain. In order to obtain a Markovian model, it is required to keep track of the type of a departing customer together with the number of customers in the system immediately after the departure of that customer. As a consequence, $(X_n,J_{n+1})$ forms a Markov chain.  \\

Taking generating functions and exploiting the fact that $X_{n-1}$ and ($A_n,J_{n+1}$) are conditionally independent, given $J_n$ and $X_{n-1}\geq 1$, we find:
\begin{align}
&\E[z^{X_{n}}1_{\{J_{n+1}=j\}}]
=\sum_{i=1}^{N}\E[z^{X_{n-1}-1}|J_{n}=i,X_{n-1} \geq 1]\E[z^{A_{n}}1_{\{J_{n+1}=j\}}|J_{n}=i,X_{n-1} \geq 1]\P(X_{n-1} \geq 1,J_{n}=i)\nonumber\\
& +\frac{B(z)}{z}\sum_{i=1}^{N}\E[z^{A_{n}}1_{\{J_{n+1}=j\}}|J_{n}=i,X_{n-1}=0]\P(X_{n-1}=0,J_{n}=i) \nonumber\\
=&\frac{1}{z}\sum_{i=1}^{N} \E[z^{X_{n-1}}1_{\{J_{n}=i\}}]\E[z^{A_{n}}1_{\{J_{n+1}=j\}}|J_{n}=i,X_{n-1} \geq 1]\nonumber\\
\ \ \ & +\frac{1}{z}\sum_{i=1}^{N}\Big(B(z)\E[z^{A_{n}}1_{\{J_{n+1}=j\}}|J_{n}=i,X_{n-1}=0]-\E[z^{A_{n}}1_{\{J_{n+1}=j\}}|J_{n}=i,X_{n-1}\geq 1]\Big)\P(X_{n-1}=0,J_{n}=i),\nonumber\\
\ \ \ &  \text{for} ~~ n=1,2,3,\dots, ~~~ j=1,2,\dots,N. \nonumber\\
\label{twentysixAA}
\end{align}
Now, we restrict ourselves to the stationary situation, assuming that the stability condition  holds.\\
Introduce, for $i,j=1,2,\dots,N$ and $|z| \leq 1$:
\begin{equation}
f_i(z)= {\rm lim}_{n \rightarrow \infty} \E[z^{X_n}1_{\{J_{n+1}=i\}}],
\label{fizA}
\end{equation}
with, for $i=1,2,\dots,N$,
\begin{equation}
f_i(0)= {\rm lim}_{n \rightarrow \infty} \P(X_n=0,J_{n+1}=i)
\label{fi0A}
\end{equation} such that
\begin{equation}
F(z)=\sum_{i=1}^{N}f_i(z). \label{F(z)}
\end{equation}
In stationarity, Equation \eqref{twentysixAA} leads to the following $N$ equations:
\begin{align}
& (z-A_{jj}(z))f_j(z)-\sum_{i=1,i\neq j}^{N}A_{ij}(z)f_i(z)
=\sum_{i=1}^{N}(B(z)A^{*}_{ij}(z)-A_{ij}(z))
f_i(0),\quad \quad ~~~~~~ j=1,2,\dots,N. \label{twentysevenA}
\end{align}
We can also write these $N$ linear equations in matrix form as
\begin{align*}
M(z)^{T}f(z)=b(z),
\end{align*}
where \begin{align}
M(z)=&
\begin{bmatrix}
z-A_{11}(z) & -A_{12}(z) & \dots & -A_{1N}(z)\\
-A_{21}(z) & z-A_{22}(z) & \dots & -A_{2N}(z)\\
\vdots & \vdots & \ddots & \vdots \\
-A_{N1}(z) & -A_{N2}(z) & \dots & z-A_{NN}(z)
\end{bmatrix}, \label{M(z)}\\
f(z)=&
\begin{bmatrix}
f_1(z)\\
f_2(z)\\
\vdots\\
f_N(z)
\end{bmatrix},
  b(z)=
\begin{bmatrix}
b_1(z)\\
b_2(z)\\
\vdots \\
b_N(z)
\end{bmatrix}, \text{ with } b_j(z)=\sum_{i=1}^{N}(B(z)A^{*}_{ij}(z)-A_{ij}(z))f_i(0).\label{b(z)}
\end{align}
Therefore, by Cramer's rule, solutions of the non-homogeneous linear system $M(z)^{T}f(z)=b(z)$ are in the form
\begin{align}
f_i(z)=\frac{\det L_i(z)}{\det M(z)^{T}}, \quad \det M(z)^{T}\neq 0,  \quad i=1,2,\dots,N, \label{f(z)_N}
\end{align} where $L_i(z)$ is the matrix formed by replacing the $i$-th column of $M(z)^T$ by the column vector $b(z)$:
\begin{align}
\det L_1(z)=&\begin{vmatrix}
         b_1(z) & -A_{21}(z)&\dots& -A_{N1}(z) \\
         b_2(z) & z-A_{22}(z)&\dots& -A_{N2}(z) \\
         \vdots & \vdots  & \ddots &\vdots \\
         b_N(z) & -A_{2N}(z)&\dots& z-A_{NN}(z) \\
       \end{vmatrix}\label{L_1(z)_D},\\
\det L_i(z)=&\begin{vmatrix}
         z-A_{11}(z)&\dots& -A_{i-11}(z)& b_1(z) & -A_{i+11}(z)&\dots& -A_{N1}(z) \\
         -A_{12}(z)&\dots& -A_{i-12}(z)& b_2(z) & -A_{i+12}(z)&\dots& -A_{N2}(z) \\
          \vdots   & \ddots &\vdots & \vdots  & \vdots & \ddots &\vdots \\
           -A_{1N}(z)&\dots& -A_{i-1N}(z)& b_N(z) & -A_{i+1N}(z)&\dots& z-A_{NN}(z) \\
       \end{vmatrix}, \quad i=2,3,\dots,N.\label{L_i(z)_D}
\end{align}

It remains to find the values of $f_1(0),f_2(0),\dots,f_N(0)$.
We shall derive $N$ linear equations for $f_1(0),f_2(0),\dots,f_N(0)$.
\\

\paragraph{First equation:}

Note that $\det M(z)^T=\det M(z)$.
After replacing the first column by sum of all $N$ columns in \eqref{M(z)}, and using \eqref{A_i(z)}, we get,
\begin{align}
\det M(z)^T=&
\begin{vmatrix}
z-A_{1}(z) & -A_{12}(z) & \dots & -A_{1N}(z)\\
z-A_{2}(z) & z-A_{22}(z) & \dots & -A_{2N}(z)\\
\vdots & \vdots & \ddots & \vdots \\
z-A_{N}(z) & -A_{N2}(z) & \dots & z-A_{NN}(z)
\end{vmatrix}
.\label{Add: detM(z)}
\end{align}
This implies that
\begin{align}\label{sumform: detM(z)}
 \det M(z)^T= \sum_{i=1}^{N}(z-A_i(z))u_{i1}(z),
\end{align}
where $u_{i1}(z)$ is the cofactor of the entry in the $i$-th row and the first column of the matrix in Equation~\eqref{Add: detM(z)}.

Note that $\{z-A_i(z)\}|_{z=1}=0, \frac{d}{dz}\{z-A_i(z)\}|_{z=1}=1-\alpha_i$, and $u_{i1}(1)=d_i$, where $d_i$ are given by Equations \eqref{d_1},\eqref{d_i},\eqref{d_N}, and $\alpha_i$ are defined in \eqref{alph_i*}, for $i=1,2,\dots,N$. Therefore, we obtain
\begin{align}\label{D: detM(z)^T}
  \frac{d}{dz}\{\det M(z)^T\}\Big|_{z=1}&=\sum_{i=1}^{N}(1-\alpha_i)d_i
  =d-\sum_{i=1}^{N}\alpha_id_i
  =d(1-\rho).
\end{align}

This implies that
\begin{align}\label{L_i(z)}
\det L_i(z)=\sum_{j=1}^Nb_j(z)r_{ji}(z), \quad i=1,2,\dots,N,
\end{align} where $b_j(z)$ is given by \eqref{b(z)}, and  $r_{ji}(z)$ is the cofactor of the entry in the $j$th row and $i$th column of the matrix $L_i(z)$, which is given by
\begin{align}
r_{11}(z)=&\begin{vmatrix}
          z-A_{22}(z)&-A_{32}(z)&\dots& -A_{N2}(z) \\
          -A_{23}(z)&z-A_{33}(z)&\dots& -A_{N3}(z) \\
          \vdots  & \vdots &\ddots &\vdots \\
          -A_{2N}(z)& -A_{3N}(z)& \dots& z-A_{NN}(z) \\
       \end{vmatrix}\label{r(11)_D},\\
 r_{1i}(z)=&(-1)^{i+1}\begin{vmatrix}
         -A_{12}(z)&\dots& -A_{i-12}(z) & -A_{i+12}(z)&\dots& -A_{N2}(z) \\
         -A_{13}(z)&\dots& -A_{i-13}(z)& -A_{i+13}(z)&\dots& -A_{N3}(z) \\
          \vdots   & \ddots &\vdots & \vdots & \ddots &\vdots \\
           -A_{1N}(z)&\dots& -A_{i-1N}(z) & -A_{i+1N}(z)&\dots& z-A_{NN}(z) \\
       \end{vmatrix}, \quad i=2,3,\dots,N,\label{r(i1)_D}\\
 r_{j1}(z)=&(-1)^{j+1}\begin{vmatrix}
         -A_{21}(z)& -A_{31}(z) & -\dots& -A_{N1}(z) \\
         \vdots   & \vdots & \ddots &\vdots \\
         -A_{2j-1}(z)& -A_{3j-1}(z)& \dots& -A_{Nj-1}(z) \\
         -A_{2j+1}(z)& -A_{3j+1}(z)& \dots& -A_{Nj+1}(z) \\
          \vdots   & \vdots & \ddots &\vdots \\
           -A_{2N}(z)& -A_{3N}(z) &\dots& z-A_{NN}(z) \\
       \end{vmatrix}, \quad j=2,3,\dots,N\label{r(1j)_D}\\
r_{ji}(z)=&(-1)^{i+j}\begin{vmatrix}
         z-A_{11}(z)&\dots& -A_{i-11}(z) & -A_{i+11}(z)&\dots& -A_{N1}(z) \\
         \vdots   & \ddots &\vdots & \vdots & \ddots &\vdots \\
         -A_{1j-1}(z)&\dots& -A_{i-1j-1}(z)& -A_{i+1j-1}(z)&\dots& -A_{Nj-1}(z) \\
         -A_{1j+1}(z)&\dots& -A_{i-1j+1}(z)& -A_{i+1j+1}(z)&\dots& -A_{Nj+1}(z) \\
          \vdots   & \ddots &\vdots & \vdots & \ddots &\vdots \\
           -A_{1N}(z)&\dots& -A_{i-1N}(z) & -A_{i+1N}(z)&\dots& z-A_{NN}(z) \\
       \end{vmatrix}, \quad i,j=2,3,\dots,N.\label{r(ij)_D}
\end{align}
Subsequently,
\begin{align}\label{L_i_prime(z)}
\frac{d}{dz}\{\det L_i(z)\}|_{z=1}&=\sum_{j=1}^N(b_j(1)r_{ji}^\prime(1)+b^\prime_j(1)r_{ji}(1))\nonumber \\
&=\sum_{j=1}^N\sum_{k=1}^{N}\Big(r_{ji}^\prime(1)(P^*_{kj}-P_{kl})+r_{ji}(1)(\E[B]P^*_{kj}+\alpha^*_{kj}-\alpha_{kj})\Big)f_k(0).
\end{align}

After replacing the first row by the sum of all $N$ rows of $\det L_i(z)$ in \eqref{L_i(z)_D}, we obtain $\det L_i(z)$, $i=2,3,\dots,N$, as
\begin{align}
\det L_i(z)=&\begin{vmatrix}
         z-A_{1}(z)&\dots& z-A_{i-1}(z)& \sum_{j=1}^N b_j(z) & z-A_{i+1}(z)&\dots& z-A_{N}(z) \\
         -A_{12}(z)&\dots& -A_{i-12}(z)& b_2(z) & -A_{i+12}(z)&\dots& -A_{N2}(z) \\
          \vdots   & \ddots &\vdots & \vdots  & \vdots & \ddots &\vdots \\
           -A_{1N}(z)&\dots& -A_{i-1N}(z)& b_N(z) & -A_{i+1N}(z)&\dots& z-A_{NN}(z) \\
       \end{vmatrix}.\label{L_i_add(z)}
\end{align}
In particular,
\begin{align*}
\det L_i(1)=&\begin{vmatrix}
         0 &\dots & 0 & 0 & 0 & \dots & 0 \\
         -P_{12}& \dots & -P_{i-12}& b_2(1) & -P_{i+12}&\dots& -P_{N2} \\
          \vdots   & \ddots &\vdots & \vdots  & \vdots & \ddots &\vdots \\
           -P_{1N}&\dots& -P_{i-1N}& b_N(1) & -P_{i+1N}&\dots& 1-P_{NN} \\
       \end{vmatrix},\quad i=2,3,\dots,N,\\
       =&\ 0.
\end{align*}
Following the same steps, one can show that  $\det L_1(1)=0$ and $\det M(1)^T=0$. Therefore, for $i=1,2,\dots,N$, we obtain,
\begin{align}
   f_i(1)&=\lim_{z\to 1}\frac{\det L_i(z)}{\det M(z)^T}\nonumber\\
   &=\frac{\frac{d}{dz}\{\det L_i(z)\}|_{z=1}}{\frac{d}{dz}\{\det M(z)^T\}|_{z=1}}.\label{f_i(1)}
\end{align}
Note that $F(1)=1$, which implies that $\sum_{i=1}^Nf_i(1)=1$. And, as a consequence, we obtain,
 \begin{align}
 \frac{\sum_{i=1}^N\frac{d}{dz}\{\det L_i(z)\}|_{z=1}}{\frac{d}{dz}\{\det M(z)^T\}|_{z=1}}&=1. \nonumber
 \end{align}
This implies that
\begin{align}
\sum_{k=1}^{N}\Bigg(\sum_{i=1}^{N}\sum_{j=1}^N\Big(r_{ji}^\prime(1)(P^*_{kj}-P_{kl})+r_{ji}(1)(\E[B]P^*_{kj}+\alpha^*_{kj}-\alpha_{kj})\Big)\Bigg)f_k(0)&=d(1-\rho).\label{first: normalizingEq}
\end{align}

\paragraph{$(N-1)$ equations:}

Under the stability condition, $\det M(z)^T$ has exactly $N-1$ zeros in $|z|<1$, denoted by $\hat{z_l}$, $l=1,2,\dots,N-1$ (see in \cite{QUESTA2017}), and $F(z)$ is an analytical function in $|z|<1$. Therefore, the numerator of $F(z)$  also has $(N-1)$ zeros in $|z|<1$. As a consequence, these $(N-1)$ zeros provide $(N-1)$ linear equations for $f_1(0),f_2(0),\dots,f_N(0)$:
\begin{align}\label{N-1: equations}
\sum_{i=1}^{N}\det L_i(\hat{z_l})&=0, \quad |\hat{z_l}|<1\nonumber\\
\implies \sum_{i=1}^{N}\sum_{j=1}^Nb_j(\hat{z_l})r_{ji}(\hat{z_l})&=0 \nonumber\\
\implies \sum_{k=1}^N\Bigg(\sum_{i=1}^{N}\sum_{j=1}^Nr_{ji}(\hat{z_l})\Big( B(\hat{z_l})A^*_{kj}(\hat{z_l})-A_{kj}(\hat{z_l})\Big)\Bigg)f_k(0)&=0,\quad l=1,2,\dots,N-1.
\end{align}

\subsection{Special case: $N=2$}

For $N=2$, we can solve \eqref{twentysevenA} and find an explicit expression for the steady-state probability generating function of the number of customers.
\begin{align}
f_1(z)=\frac{\sum_{i=1}^{2}\Big((z-A_{22}(z))(B(z)A^{*}_{i1}(z)-A_{i1}(z))+A_{21}(z)(B(z)A^{*}_{i2}(z)-A_{i2}(z)) \Big)f_i(0)}{(z-A_{11}(z))(z-A_{22}(z))-A_{12}(z)A_{21}(z)}, \label{twentynine}
\end{align}
\begin{align}
f_2(z)=\frac{\sum_{i=1}^{2}\Big((z-A_{11}(z))(B(z)A^{*}_{i2}(z)-A_{i2}(z))+A_{12}(z)(B(z)A^{*}_{i1}(z)-A_{i1}(z)) \Big)f_i(0)}{(z-A_{11}(z))(z-A_{22}(z))-A_{12}(z)A_{21}(z)}. \label{thirty}
\end{align}
In particular,
\begin{align}
&f_1(1)=\lim_{z\to 1}\frac{\sum_{i=1}^{2}\Big((z-A_{22}(z))(B(z)A^{*}_{i1}(z)-A_{i1}(z))+A_{21}(z)(B(z)A^{*}_{i2}(z)-A_{i2}(z)) \Big)f_i(0)}{(z-A_{11}(z))(z-A_{22}(z))-A_{12}(z)A_{21}(z)}\nonumber\\
=&\frac{\sum_{i=1}^{2}\Big(P_{21}(P^{*}_{i1}\E[B]+\alpha^{*}_{i1}-\alpha_{i1}) +(1-\alpha_{22})(P^*_{i1}-P_{i1}) +P_{21}(P^{*}_{i2}\E[B]+\alpha^{*}_{i2}-\alpha_{i2}) +\alpha_{21}(P^*_{i2}-P_{i2}) \Big)f_i(0)}{(1-P_{11})(1-\alpha_{22})+(1-P_{22})(1-\alpha_{11})-P_{12}\alpha_{21}-P_{21}\alpha_{12}}\nonumber\\
=&\frac{\sum_{i=1}^{2}\Big(P_{21}(\E[B]+\alpha^*_i-\alpha_i)+(1-\alpha_2)(P^*_{i1}-P_{i1})\Big)f_i(0)}{(P_{12}+P_{21})\left(1-\frac{P_{21}}{P_{12}+P_{21}}\alpha_1-\frac{P_{12}}{P_{12}+P_{21}}\alpha_2\right)}\nonumber\\
=&\frac{\sum_{i=1}^{2}\Big(P_{21}(\E[B]+\alpha^*_i-\alpha_i)+(1-\alpha_2)(P^*_{i1}-P_{i1})\Big)f_i(0)}{(P_{12}+P_{21})\left(1-\rho \right)}\label{f_1(1)}.
\end{align}
Similarly,
 \begin{align}
f_2(1)=\frac{\sum_{i=1}^{2}\Big(P_{12}(\E[B]+\alpha^*_i-\alpha_i)+(1-\alpha_1)(P^*_{i2}-P_{i2})\Big)f_i(0)}{(P_{12}+P_{21})\left(1-\rho \right)}\label{f_2(1)}.
\end{align}
As a consequence of $f_1(1)+f_2(1)=1$, we obtain
\begin{align}\label{F(1)}
  \sum_{i=1}^{2}\Big((P_{12}+P_{21})(\E[B]+\alpha^*_i-\alpha_i)+(\alpha_1-\alpha_2)(P^*_{i1}-P_{i1})\Big)f_i(0)=(P_{12}+P_{21})\left(1-\rho \right).
\end{align}
After substituting the values of $f_1(z)$ and $f_2(z)$ from Equations \eqref{twentynine} and \eqref{thirty}, respectively, in \eqref{F(z)}, we obtain
\begin{align}
F(z)=\frac{\sum_{i=1}^{2}\Big((z+A_{12}(z)-A_{22}(z))(B(z)A^{*}_{i1}(z)-A_{i1}(z))+(z+A_{21}(z)-A_{11}(z))(B(z)A^{*}_{i2}(z)-A_{i2}(z)) \Big)f_i(0)}{(z-A_{11}(z))(z-A_{22}(z))-A_{12}(z)A_{21}(z)}.\label{F_2(z)}
\end{align}
Let $z=\hat{z}$ be the zero of the denominator of $F(z)$ such that $|\hat{z}|<1$. Since $z=\hat{z}$ must also be the zero of the numerator of $F(z)$, we obtain the following equation in terms of $f_1(0)$ and $f_2(0)$:
\begin{align}\label{second: eq}
\sum_{i=1}^{2}\Big((\hat{z}+A_{12}(\hat{z})-A_{22}(\hat{z}))(B(\hat{z})A^{*}_{i1}(\hat{z})-A_{i1}(\hat{z}))+(\hat{z}+A_{21}(\hat{z})-A_{11}(\hat{z}))(B(\hat{z})A^{*}_{i2}(\hat{z})-A_{i2}(\hat{z})) \Big)f_i(0)=0.
\end{align}
Solving Equations \eqref{F(1)}and \eqref{second: eq} yields
\begin{align}
  f_1(0)=&\frac{-(P_{12}+P_{21})\left(1-\rho \right)R_{12}}{\det R}, \label{N: f1(0)}\\
  f_2(0)=&\frac{(P_{12}+P_{21})\left(1-\rho \right)R_{11}}{\det R}, \label{N: f2(0)}
\end{align} where $\det R$ is the determinant of the matrix $R=[R_{ij}]$, whose elements are given by
\begin{align*}
R_{1j}&=(\hat{z}+A_{12}(\hat{z})-A_{22}(\hat{z}))(B(\hat{z})A^{*}_{j1}(\hat{z})-A_{j1}(\hat{z}))+(\hat{z}+A_{21}(\hat{z})-A_{11}(\hat{z}))(B(\hat{z})A^{*}_{j2}(\hat{z})-A_{j2}(\hat{z})),\\
R_{2j}&=(P_{12}+P_{21})(\E[B]+\alpha^*_j-\alpha_j)+(\alpha_1-\alpha_2)(P^*_{j1}-P_{j1}), \quad j=1,2.
\end{align*}

\subsection{Stationary queue length analysis: arrival and arbitrary epochs}

In the previous subsection, we determined the probability generating function of the stationary queue length distribution at the departure epoch of an arbitrary customer for general batch arrivals. As customers arrive at the system according to a batch Poisson process with rate $\lambda$, from the PASTA property, the distribution of the number of customers in the system at the arrival time of a batch is the same as the distribution of the number of customers at an arbitrary time. After using PASTA and level-crossing arguments (see \cite{QUESTA2017} for more details), we obtain the following relations:
\begin{align}
 \E[z^{X}]=\E[z^{X^{\textit{ca}}}] = \E[z^{X^{\textit{ba}}}]\frac{1-B(z)}{\E[B](1-z)}, \label{rel:queue}
  \end {align}
  with,
 \begin{align}
 \E[z^{X^{\textit{arb}}}]=\E[z^{X^{\textit{ba}}}] ,
 \label{F(z):arbitrary}
\end{align}
where $X$ and $X^{\textit{ca}}$ are the number of customers at the departure and the arrival epoch of the customer respectively;   $X^{\textit{arb}}$ and $X^{\textit{ba}}$ are the number of customers at an arbitrary time and the arrival time of a batch respectively.
From these relations, we can obtain all the required distributions.

\section{Heavy-traffic analysis}\label{heavy-traffc}

In this section, we shall determine the HT limit of the scaled queue length at departure epochs. In particular, we will show that under some conditions the distribution of the scaled stationary queue length in heavy traffic is exponential. This will be formally stated in Theorem \ref{thm:ch6}.

Let us define the HT limit for the LST of the scaled queue length at departure epochs, $(1-\rho)X$, for $i=1,2,\dots,N$:
\begin{align*}
\bar{F}(s)=\lim_{\rho \uparrow 1}\E[e^{-s(1-\rho)X}]=\sum_{i=1}^N\bar{f}_i(s),
\end{align*}
with
\begin{align*}
\bar{f}_i(s)=\lim_{\rho \uparrow 1}\E[e^{-s(1-\rho)X_n}1_{\{J_{n+1}=i\}}].
\end{align*}
Firstly, we introduce the following notation. For $i=1,2,\dots,N$, $|z|\leq 1$,
\begin{align}\label{mathcal: A_i}
\mathcal{A}_{i}(z)=\E[z^{A_n}|J_{n+1}=i,X_{n-1}\geq 1],
\end{align}
and
\begin{align}\label{gamma_i}
\gamma_i=\E[A_n|J_{n+1}=i,X_{n-1}\geq 1].
\end{align}
To find the limiting HT distribution, we substitute
\begin{align}
z=e^{-s(1-\rho)}=1-s(1-\rho)+\frac{1}{2}s^2(1-\rho)^2+O((1-\rho)^3), ~~~~~\text{as } \rho \uparrow 1.\label{HT_z}
\end{align}
With this substitution, we can write the generating function of the batch size, $B(z)$, as
\begin{align}
B(e^{-s(1-\rho)})=1-s(1-\rho)\E[B]+\frac{1}{2}s^2(1-\rho)^2\E[B^2]+O((1-\rho)^3). \label{HT_B(z)}
\end{align}
Similarly, from Equations \eqref{A_ij(z)}, \eqref{A_ijstar(z)}, and \eqref{mathcal: A_i}, we obtain for $\rho \uparrow 1$,
\begin{align}
A_{ij}(e^{-s(1-\rho)})=&P_{ij}-s(1-\rho)\alpha_{ij}+\frac{1}{2}s^2(1-\rho)^2\hat{\alpha}_{ij}+O((1-\rho)^3),\label{HT_A_ij}\\
A^*_{ij}(e^{-s(1-\rho)})=&P^*_{ij}-s(1-\rho)\alpha^*_{ij}+\frac{1}{2}s^2(1-\rho)^2\hat{\alpha}^*_{ij}+O((1-\rho)^3), \label{HT_A^*_ij}\\
\mathcal{A}_{i}(e^{-s(1-\rho)})=&1-s(1-\rho)\gamma_{i}+\frac{1}{2}s^2(1-\rho)^2\hat{\gamma}_{i}+O((1-\rho)^3),~\text{as } \rho \uparrow 1,\label{HT_mathcal_A_i}
\end{align} where $\alpha_{ij},\alpha^*_{ij}$, and $\gamma_i$ are respectively defined in Equations \eqref{alpha_ij}, \eqref{alpha_ij*} and \eqref{gamma_i}, and
\begin{align}
\hat{\alpha}_{ij}=&\E[(A_n)^2 1_{\{J_{n+1}=j\}}|J_n=i,X_{n-1}\geq 1],\label{hat: alpha_ij}\\
\hat{\alpha}^*_{ij}=&\E[(A^*_n)^2 1_{\{J_{n+1}=j\}}|J_n=i,X_{n-1}=0],\label{hat: alpha^*_ij}\\
\hat{\gamma}_i=&\E[(A_n)^2|J_{n+1}=i,X_{n-1}\geq 1].\label{hat: gamma_i}
\end{align}
Using $A_i(z)=\sum_{j=1}^{N}A_{ij}(z)$ and $A^*_i(z)=\sum_{j=1}^{N}A^*_{ij}(z)$, we obtain,
\begin{align}
A_{i}(e^{-s(1-\rho)})=&1-s(1-\rho)\alpha_{i}+\frac{1}{2}s^2(1-\rho)^2\hat{\alpha}_{i}+O((1-\rho)^3),\label{HT_A_i}\\
A^*_{i}(e^{-s(1-\rho)})=&1-s(1-\rho)\alpha^*_{i}+\frac{1}{2}s^2(1-\rho)^2\hat{\alpha}^*_{i}+O((1-\rho)^3),~\text{as } \rho \uparrow 1, \label{HT_A^*_i}
\end{align}
where $\alpha_{i}$ and $\alpha^*_{i}$ are respectively defined in Equations \eqref{alph_i} and \eqref{alph_i*}, and
\begin{align}
\hat{\alpha}_i=&\sum_{j=1}^N\hat{\alpha}_{ij},\label{hat: alpha_i}\\
\hat{\alpha}^{*}_i=&\sum_{j=1}^N\hat{\alpha}^{*}_{ij}.\label{hat: alpha^*_i}
\end{align}

After substituting the values of $B(z),A_{ij}(z)$ and $A^*_{ij}(z)$ from Equations \eqref{HT_B(z)}, \eqref{HT_A_ij} and \eqref{HT_A^*_ij} respectively, we obtain $b_j(z)$, with $z=e^{-s(1-\rho)}$, from  \eqref{b(z)}  as
\begin{align}
&b_j(e^{-s(1-\rho)})=\sum_{i=1}^N\Big(\ (P^*_{ij}-P_{ij})-s(1-\rho)(P^*_{ij}\E[B]+\alpha^*_{ij}-\alpha_{ij})\nonumber\\
&+\frac{s^2(1-\rho)^2}{2}(\E[B^2]P^*_{ij}+2\E[B]\alpha^*_{ij}+\hat{\alpha}^*_{ij}-\hat{\alpha}_{ij})+O((1-\rho)^3)\Big)f_i(0),~\text{as } \rho \uparrow 1.\label{HT_b_j}
\end{align}
Summing over $j$ and using $\sum_{j=1}^NP_{ij}=\sum_{j=1}^N P^*_{ij}=1$ gives
\begin{align}
\sum_{j=1}^N b_j(e^{-s(1-\rho)})&=-s(1-\rho)\sum_{i=1}^N\Big((\E[B]+\alpha^*_{i}-\alpha_{i})-\frac{s(1-\rho)}{2}(\E[B^2]\nonumber\\
&+2\E[B]\alpha^*_{i}+\hat{\alpha}^*_{i}-\hat{\alpha}_{i})+O((1-\rho)^2)\Big)f_i(0), ~\text{as } \rho \uparrow 1.\label{HT_sum_b_j}
\end{align}

Substituting the values of $z,A_{ij}(z),A_{i}(z),b_j(z)$, and $\sum_{j=1}^{N}b_j(z)$ from Equations \eqref{HT_z}, \eqref{HT_A_ij}, \eqref{HT_A_i}, \eqref{HT_b_j}, and \eqref{HT_sum_b_j}, respectively, in Equation \eqref{L_i_add(z)}, and after simplification, we obtain $\det L_i(z)$, with $z=e^{-s(1-\rho)}$, as
\begin{align}
&\det L_i(e^{-s(1-\rho)})=-s(1-\rho)\nonumber\\
&\times\begin{vmatrix}
         1-\alpha_{1}&\dots& 1-\alpha_{i-1}& \sum_{k=1}^N(\E[B]+\alpha^*_{k}-\alpha_{k})f_k(0)  & 1-\alpha_{i+1}&\dots&1-\alpha_{N}\\
         -P_{12}&\dots& -P_{i-12}& \sum_{k=1}^N(P^*_{k2}-P_{k2})f_k(0)  & -P_{i+12}&\dots& -P_{N2} \\
          \vdots   & \ddots &\vdots & \vdots  & \vdots & \ddots &\vdots \\
           -P_{1N}&\dots& -P_{i-1N}&\sum_{k=1}^N(P^*_{kN}-P_{kN})f_k(0) & -P_{i+1N}&\dots& 1-P_{NN} \\
       \end{vmatrix}\nonumber\\
       &+c_is^2(1-\rho)^2+O((1-\rho)^3),\quad i=2,3,\dots,N,\label{HT_Li1}
\end{align} where $c_i$ is the coefficient of $s^2(1-\rho)^2$ term such that
\begin{align}\label{lim_ci}
\lim_{\rho \uparrow 1} c_i=0.
\end{align}
This coefficient exists, because
\[\lim_{\rho \uparrow 1}\P(X=0)= \lim_{\rho \uparrow 1} \sum_{k=1}^{N} f_k(0) =0,\]
implying that  $\lim_{\rho \uparrow 1} f_k(0)=0$ for all $1\leq k  \leq N$.

Now, differentiating Equation \eqref{L_i_add(z)} w.r.t. $z$, and substituting $z=1$, we get,
 \begin{align}\label{HT_Li2}
 &\frac{d}{dz}\{\det L_i(z)\}|_{z=1}=\begin{vmatrix}
         1-\alpha_{1}&\dots& 1-\alpha_{i-1}& \sum_{k=1}^N(\E[B]+\alpha^*_{k}-\alpha_{k})f_k(0)  & 1-\alpha_{i+1}&\dots&1-\alpha_{N}\\
         -P_{12}&\dots& -P_{i-12}& \sum_{k=1}^N(P^*_{k2}-P_{k2})f_k(0)  & -P_{i+12}&\dots& -P_{N2} \\
          \vdots   & \ddots &\vdots & \vdots  & \vdots & \ddots &\vdots \\
           -P_{1N}&\dots& -P_{i-1N}&\sum_{k=1}^N(P^*_{kN}-P_{kN})f_k(0) & -P_{i+1N}&\dots& 1-P_{NN} \\
       \end{vmatrix}.
 \end{align}
After using Equations \eqref{D: detM(z)^T}, \eqref{f_i(1)} and \eqref{HT_Li2} in Equation \eqref{HT_Li1}, we can write
\begin{align}
&\det L_i(e^{-s(1-\rho)})\nonumber\\
&=-s(1-\rho)\frac{d}{dz}\{\det L_i(z)\}|_{z=1}+c_is^2(1-\rho)^2+O((1-\rho)^3)\nonumber\\
&=-s(1-\rho)\frac{d}{dz}\{\det M(z)^T\}|_{z=1}f_i(1)+c_is^2(1-\rho)^2+O((1-\rho)^3)\nonumber\\
&=-sd(1-\rho)^2(f_i(1)-\frac{c_is}{d})+O((1-\rho)^3),\quad i=2,3,\dots,N.\nonumber
\end{align}
Similarly,
\begin{align*}
\det L_1(e^{-s(1-\rho)})&=-sd(1-\rho)^2(f_1(1)-\frac{c_1s}{d})+O((1-\rho)^3).
\end{align*}
Hence, we can write, for $i=1,2,\dots,N$,
\begin{align}
\det L_i(e^{-s(1-\rho)})&=-sd(1-\rho)^2(f_i(1)-\frac{c_is}{d})+O((1-\rho)^3).\label{HT_Li(s)}
\end{align}
From Equation \eqref{Add: detM(z)}, $\det M(z)^T$ is given by
\begin{align}
&\det M(z)^T=
\begin{vmatrix}
z-A_{1}(z) & -A_{12}(z) & \dots & -A_{1N}(z)\\
z-A_{2}(z) & z-A_{22}(z) & \dots & -A_{2N}(z)\\
\vdots & \vdots & \ddots & \vdots \\
z-A_{N}(z) & -A_{N2}(z) & \dots & z-A_{NN}(z)
\end{vmatrix}\nonumber\\
&=\frac{1}{\prod_{i=1}^{N}\pi_i}
\begin{vmatrix}
\pi_1(z-A_{1}(z)) & -\pi_1 A_{12}(z) & \dots & -\pi_1 A_{1N}(z)\\
\pi_2(z-A_{2}(z)) & \pi_2(z-A_{22}(z)) & \dots & -\pi_2A_{2N}(z)\\
\vdots & \vdots & \ddots & \vdots \\
\pi_N(z-A_{N}(z)) & -\pi_NA_{N2}(z) & \dots & \pi_N(z-A_{NN}(z))
\end{vmatrix},\nonumber\\
 &\qquad \qquad \qquad \qquad \qquad\text{ since } \pi_i\neq 0, i=1,2,\dots,N.\label{HT: M(z)0}
\end{align}
Using $\lim_{\rho \uparrow 1}f_k(0)=0$,  we will first show that $\lim_{\rho \uparrow 1}f_j(1)=\pi_j$ for all $1\leq j\leq N$. To do so, we first write $\lim_{\rho \uparrow 1}f_j(1)$ as
\begin{align*}
\lim_{\rho \uparrow 1}f_j(1)&=\lim_{\rho \uparrow 1}\P(J_{n+1}=j)\\
&=\lim_{\rho \uparrow 1}(\P(J_{n+1}=j,X_{n-1}=0)+\P(J_{n+1}=j,X_{n-1}\geq 1))\\
&=\lim_{\rho \uparrow 1}\P(J_{n+1}=j|X_{n-1}\geq 1)\P(X_{n-1}\geq 1)\\
&=\lim_{\rho \uparrow 1}\sum_{i=1}^{N}\P(J_{n+1}=j|J_n=i,X_{n-1}\geq 1)\P(J_{n}=i|X_{n-1}\geq 1)\P(X_{n-1}\geq 1)\\
&=\lim_{\rho \uparrow 1}\sum_{i=1}^{N}P_{ij}P(J_{n}=i,X_{n-1}\geq 1)\\
&=\lim_{\rho \uparrow 1}\sum_{i=1}^{N}P_{ij}(\P(J_n=i)-\P(J_{n}=i,X_{n-1}=0))\\
&=\lim_{\rho \uparrow 1}\sum_{i=1}^{N}P_{ij}(f_i(1)-f_i(0))\\
&=\sum_{i=1}^{N}P_{ij}\lim_{\rho \uparrow 1}f_i(1), \quad \text{for } j=1,2,\dots,N.
\end{align*}
As $P=[P_{ij}]_{i,j\in\{1,2,\dots,N\}}$ is the transition probability matrix of an irreducible discrete time Markov chain, with stationary distribution $\pi=(\pi_1,\pi_2,\dots,\pi_N)$, $\pi$ is the unique solution of the system of equations  $\pi (I-P)=0$, and,  hence, $\lim_{\rho \uparrow 1}f_j(1)=\pi_j$ for all $1\leq j\leq N$. As a consequence, we obtain $\lim_{\rho \uparrow 1} \P(J_{n}=j|X_{n-1}\geq 1)=\pi_j$.

Furthermore,
\begin{align}
  \lim_{\rho \uparrow 1}\P(J_{n+1}=j|X_{n-1}\geq 1)&=\lim_{\rho \uparrow 1}\sum_{i=1}^{N}\P(J_{n+1}=j|J_n=i,X_{n-1}\geq 1)\P(J_{n}=i|X_{n-1}\geq 1)\nonumber\\
  &=\sum_{i=1}^{N}P_{ij}\pi_i \nonumber\\
  &=\pi_j.\label{P(J_{n+1}|X>0)}
\end{align}
As a consequence, $\lim_{\rho \uparrow 1} \mathcal{A}_j(z)$ is given by
\begin{align}
\lim_{\rho \uparrow 1} \mathcal{A}_j(z)&=\lim_{\rho \uparrow 1}\frac{\E[z^{A_n}1_{\{J_{n+1}=j\}}|X_{n-1}\geq 1]}{\P(J_{n+1}=j|X_{n-1}\geq 1)}\nonumber\\
&=\lim_{\rho \uparrow 1}\frac{\sum_{i=1}^{N} \P(J_{n}=i|X_{n-1}\geq 1) \E[z^{A_n}1_{\{J_{n+1}=j\}}|J_n=i,X_{n-1}\geq 1]}{\P(J_{n+1}=j|X_{n-1}\geq 1)} \nonumber\\
&=\frac{\sum_{i=1}^{N} \pi_i A_{ij}(z)}{\pi_j}. \label{Rel: mathcalA_jA_ij}
\end{align}
Subsequently, we obtain,
\begin{align}
\gamma_j=\frac{\sum_{i=1}^{N} \pi_i \alpha_{ij}}{\pi_j},\quad \text{as } \rho \uparrow 1.\label{Rel: gammaj_alphaij}
\end{align}
Replacing the first row by the sum of all $N$ rows in Equation \eqref{HT: M(z)0}, and using $\mathcal{A}_j(z)=\frac{1}{\pi_j}\sum_{i=1}^{N}\pi_iA_{ij}(z)$ as $\rho \uparrow 1$ and $\sum_{i=1}^{N}\pi_i=1$ , we obtain $\det M(z)^T$ as, for $\rho \uparrow 1$,
\begin{align}
\det M(z)^T&=\frac{1}{\prod_{i=1}^{N}\pi_i}
\begin{vmatrix}
z-\sum_{i=1}^{N}\pi_iA_i(z) & \pi_2(z- \mathcal{A}_{2}(z)) & \dots & \pi_N(z-\mathcal{A}_{N}(z))\\
\pi_2(z-A_{2}(z)) & \pi_2(z-A_{22}(z)) & \dots & -\pi_2A_{2N}(z)\\
\vdots & \vdots & \ddots & \vdots \\
\pi_N(z-A_{N}(z)) & -\pi_NA_{N2}(z) & \dots & \pi_N(z-A_{NN}(z))
\end{vmatrix}\nonumber\\
&=\frac{1}{\pi_1}
\begin{vmatrix}
z-\sum_{i=1}^{N}\pi_iA_i(z) & \pi_2(z-\mathcal{A}_{2}(z)) & \dots & \pi_N(z-\mathcal{A}_{N}(z))\\
z-A_{2}(z) & z-A_{22}(z) & \dots & -A_{2N}(z)\\
\vdots & \vdots & \ddots & \vdots \\
z-A_{N}(z) & -A_{N2}(z) & \dots & z-A_{NN}(z)
\end{vmatrix}.\label{HT_M(z)1}
\end{align}
Substituting the values of $z,A_{ij}(z),A_{i}(z)$, and $\mathcal{A}_{i}(z)$ from Equations \eqref{HT_z}, \eqref{HT_A_ij}, \eqref{HT_A_i}, and \eqref{HT_mathcal_A_i}, respectively, in Equation \eqref{HT_M(z)1}, and after simplification, with $z=e^{-s(1-\rho)},\rho=\sum_{i=1}^{N}\pi_i\alpha_i,\hat{\alpha}=\sum_{i=1}^{N}\pi_i\hat{\alpha}_i,\pi_1=\frac{d_1}{d}$, we obtain
\begin{align}
&\det M(e^{-s(1-\rho)})^T\nonumber \\
&=\frac{d}{d_1}
\begin{vmatrix}
-s(1-\rho)^2(1-\frac{s}{2}(1-\hat{\alpha})) & -\pi_2s(1-\rho)(1-\gamma_2)& \dots & -\pi_Ns(1-\rho)(1-\gamma_N)\\
-s(1-\rho)(1-\alpha_2) & 1-P_{22} & \ddots & -P_{2N}\\
\vdots & \vdots & \ddots & \vdots \\
-s(1-\rho)(1-\alpha_N) & -P_{N2} & \ddots & 1-P_{NN}
\end{vmatrix}+O((1-\rho)^3)\nonumber\\
&=\frac{-sd(1-\rho)^2}{d_1} \begin{vmatrix}
1-\frac{s}{2}(1-\hat{\alpha}) & -\pi_2s(1-\gamma_2)& \dots & -\pi_Ns(1-\gamma_N)\\
1-\alpha_2 & 1-P_{22} & \ddots & -P_{2N}\\
\vdots & \vdots & \ddots & \vdots \\
1-\alpha_N & -P_{N2} & \ddots & 1-P_{NN}
\end{vmatrix}
+O((1-\rho)^3)\nonumber\\
&=\frac{-sd(1-\rho)^2}{d_1}\Big((1-\frac{s}{2}(1-\hat{\alpha}))d_1-s\sum_{k=2}^{N}\pi_k(1-\gamma_k)q_k\Big)+O((1-\rho)^3)\nonumber\\
&=-sd(1-\rho)^2\left(1+s\left(\frac{\hat{\alpha}-1}{2}-\frac{1}{d_1}\sum_{k=2}^{N}\pi_k(1-\gamma_k)q_k\right)\right)+O((1-\rho)^3),\label{MT_M(s)^T}
\end{align} where $d_1$ is defined in Equation \eqref{d_1}, and $q_k$, $k=2,3,\dots,N$, is the cofactor of the entry in the first row and the $k$-th column of the matrix
\[\begin{bmatrix}
1-\frac{s}{2}(1-\hat{\alpha}) & -\pi_2s(1-\gamma_2)& \dots & -\pi_Ns(1-\gamma_N)\\
1-\alpha_2 & 1-P_{22} & \ddots & -P_{2N}\\
\vdots & \vdots & \ddots & \vdots \\
1-\alpha_N & -P_{N2} & \ddots & 1-P_{NN}
\end{bmatrix},\]
which is given by
\begin{align}
q_2=&-\begin{vmatrix}
1-\alpha_2 & -P_{23} & -P_{24} &  \dots &-P_{2N}\\
1-\alpha_3 & 1-P_{33} &  -P_{34} & \dots &1-P_{3N}\\
\vdots & \vdots & \vdots  &\ddots & \vdots\\
1-\alpha_N & -P_{N3} & -P_{N4} & \dots &1-P_{NN}
\end{vmatrix},\label{q1}\\
q_k=&(-1)^{k+1}\begin{vmatrix}
1-\alpha_2 & 1-P_{22} & \dots & -P_{2k-1} & -P_{2K+1}& \dots &-P_{2N}\\
1-\alpha_3 & -P_{32} & \dots & -P_{3k-1} & -P_{3K+1}& \dots &1-P_{3N}\\
\vdots & \vdots & \ddots & \vdots  & \vdots &\ddots & \vdots\\
1-\alpha_N & -P_{N2} & \dots & -P_{Nk-1} & -P_{NK+1}& \dots &1-P_{NN}
\end{vmatrix},\label{qk}
\end{align}
for $k=3,4,\dots,N$.
Therefore, 
\begin{align}
\bar{f}_i(s)&=\lim_{\rho \uparrow 1} \frac{\det L_i(e^{-s(1-\rho)})}{\det M(e^{-s(1-\rho)})^T}\nonumber\\
&=\lim_{\rho \uparrow 1} \frac{-sd(1-\rho)^2(f_i(1)-\frac{c_is}{d})+O((1-\rho)^3)}{-sd(1-\rho)^2\left(1+s\left(\frac{\hat{\alpha}-1}{2}-\frac{1}{d_1}\sum_{k=2}^{N}\pi_k(1-\gamma_k)q_k\right)\right)+O((1-\rho)^3)}\nonumber\\
&=\frac{\pi_i}{1+s\left(\frac{\bar{\hat{\alpha}}-1}{2}-\frac{1}{d_1}\sum_{k=2}^{N}\pi_k(1-\bar{\gamma}_k)\bar{q}_k\right)},\label{HT_fi(s)}
\end{align} where we define $\lim_{\rho \uparrow 1}\hat{\alpha}=\bar{\hat{\alpha}},\lim_{\rho \uparrow 1}\gamma_k=\bar{\gamma}_k$ and $\lim_{\rho \uparrow 1}q_k=\bar{q}_k$.\\


This finally gives us the HT limit of the scaled queue length, which we formulate in the following theorem.

\begin{theorem}\label{thm:ch6}
If $\E[B^2]$ and $\hat{\alpha}_{ij}$ are finite for $i,j=1,2,\dots,N$, then
\begin{align}
\bar{F}(s)= \lim_{\rho \uparrow 1} \E[e^{-s(1-\rho)X}]= \frac{1}{1+s\left(\frac{\bar{\hat{\alpha}}-1}{2}-\frac{1}{d_1}\sum_{k=2}^{N}\pi_k(1-\bar{\gamma}_k)\bar{q}_k\right)},
\label{HT_F(s)}
\end{align}
provided $\left(\frac{\bar{\hat{\alpha}}-1}{2}-\frac{1}{d_1}\sum_{k=2}^{N}\pi_k(1-\bar{\gamma}_k)\bar{q}_k\right)>0$, which is the LST of an exponentially distributed random variable with rate parameter $$\eta=\frac{1}{\frac{\bar{\hat{\alpha}}-1}{2}-\frac{1}{d_1}\sum_{k=2}^{N}\pi_k(1-\bar{\gamma}_k)\bar{q}_k}.$$
\end{theorem}
\begin{remark}\label{remark: equal_alpha_i}
If $\alpha_i=\alpha$ for all $i=1,2,\dots,N$, then Equation \eqref{rho} implies that $\rho=\alpha$. In that case, the system is in heavy traffic when $\alpha\uparrow 1$  and, as a consequence, when $\alpha_i\uparrow 1$ for all $i=1,2,\dots,N$.
Note that each element of the first column of $q_k$, $k=2,3,\dots,N$, tends to zero, as $\alpha_i\uparrow 1$ for all $i=1,2,\dots,N$. It follows that $q_k=0$, which implies that $\bar{q}_k=0$  for all $k=2,3,\dots,N$, and
\begin{align*}
\bar{F}(s) \rightarrow \frac{1}{1+s\left(\frac{\bar{\hat{\alpha}}-1}{2}\right)},\quad \text{as }\alpha=\rho \uparrow 1,
\end{align*} which is the HT limit of the scaled queue length of the standard $M^X/G/1$ without dependencies at the departure epochs. Furthermore, we can conclude that the term $-\frac{s}{d_1}\sum_{k=2}^{N}\pi_k(1-\bar{\gamma}_k)\bar{q}_k$ in Equation \eqref{HT_F(s)} appears due to the dependent service times.
\end{remark}
\begin{remark}\label{remark: N=2}
For $N=2$, Equation~\eqref{HT_F(s)} reduces to
\begin{align*}
\bar{F}(s)=\frac{1}{1+s\left(\frac{\bar{\hat{\alpha}}-1}{2}+\frac{\left(1-\bar{\alpha}_2\right)}{P_{12}+P_{21}}\left(\frac{P_{12}}{P_{21}}(1-\bar{\alpha}_{22})-\bar{\alpha}_{12} \right)\right)}.
\end{align*} 
Additionally, when $\frac{\left(1-\bar{\alpha}_2\right)}{P_{12}+P_{21}}\left(\frac{P_{12}}{P_{21}}(1-\bar{\alpha}_{22})-\bar{\alpha}_{12} \right)=0$,
then $\bar{F}(s)$ becomes $(1+s((\bar{\hat{\alpha}}-1)/{2}))^{-1}$, which is the
HT limit of the scaled queue length at departure epochs of the standard $M^X/G/1$ queue without dependencies.
\end{remark}
\begin{remark}
After using Equation \eqref{HT_F(s)} in \eqref{rel:queue} and \eqref{F(z):arbitrary}, it can be shown by substituting $z=e^{-s(1-\rho)}$ and taking $\rho\uparrow 1$ that the HT distribution of the scaled stationary queue length at an \emph{arbitrary epoch} is the same as  the HT distribution of the scaled stationary queue length at a departure epoch.
\end{remark}

\section{Numerical example}\label{numerical_results}


In this section we would like to given an example of the interesting situation described in Remark~\ref{remark: N=2}, where we carefully construct the dependencies between successive service times in such a way that they disappear as $\rho$ tends to one.
For simplicity, we take $N=2$, $B(z)=z$ and $\tilde{G}^*_{ij}(s)=\tilde{G}_{ij}(s)$ for all $i,j=1,2$, i.e., there are two customer types, the batch size is one, and customers arriving in an empty system have the same service-time distributions as regular customers. The conditional service times are Erlang distributed random variables, with
\begin{align*}
G_{ij}(x) &=  \Big(1-\sum_{m=0}^{k_{ij}-1} \frac{(\mu_{ij} x)^m}{m!}e^{-\mu_{ij}x}\Big)P_{ij},
\end{align*}
where $k_{ij}=i+j,\mu_{ij} > 0,$ $i,j=1, 2$.
We can use Equation \eqref{rel_A_G} to obtain
\begin{align*}
A_{ij}(z)=P_{ij}\left(\frac{\mu_{ij}}{\lambda(1-B(z))+\mu_{ij}}\right)^{k_{ij}}, \quad \text{for } i,j=1,2.
\end{align*}
We choose model parameters $P_{11}=0.9$, $\alpha_{11}=\lambda, \alpha_{12}=3\lambda,\alpha_{21}=10 \lambda$, and $\alpha_{22}=20 \lambda$. To ensure that $\frac{P_{12}}{P_{21}}(1-\bar{\alpha}_{22})-\bar{\alpha}_{12}=0$,
we take $P_{22}=0.951138$.
\begin{figure}[ht]
\begin{center}
\includegraphics[width=\textwidth]{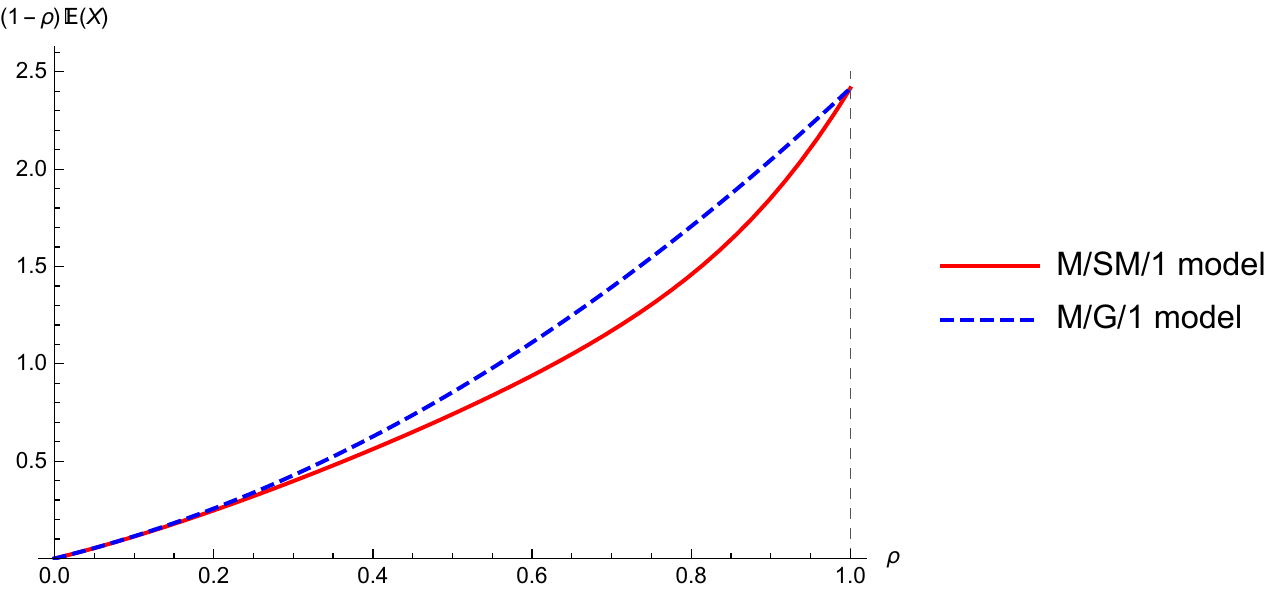}
\end{center}
\caption{The mean scaled queue length versus the number of arrivals per time unit.}
\label{fig:EXscaledexample}
\end{figure}

\begin{figure}[ht!]
\begin{center}
\includegraphics[width=\textwidth]{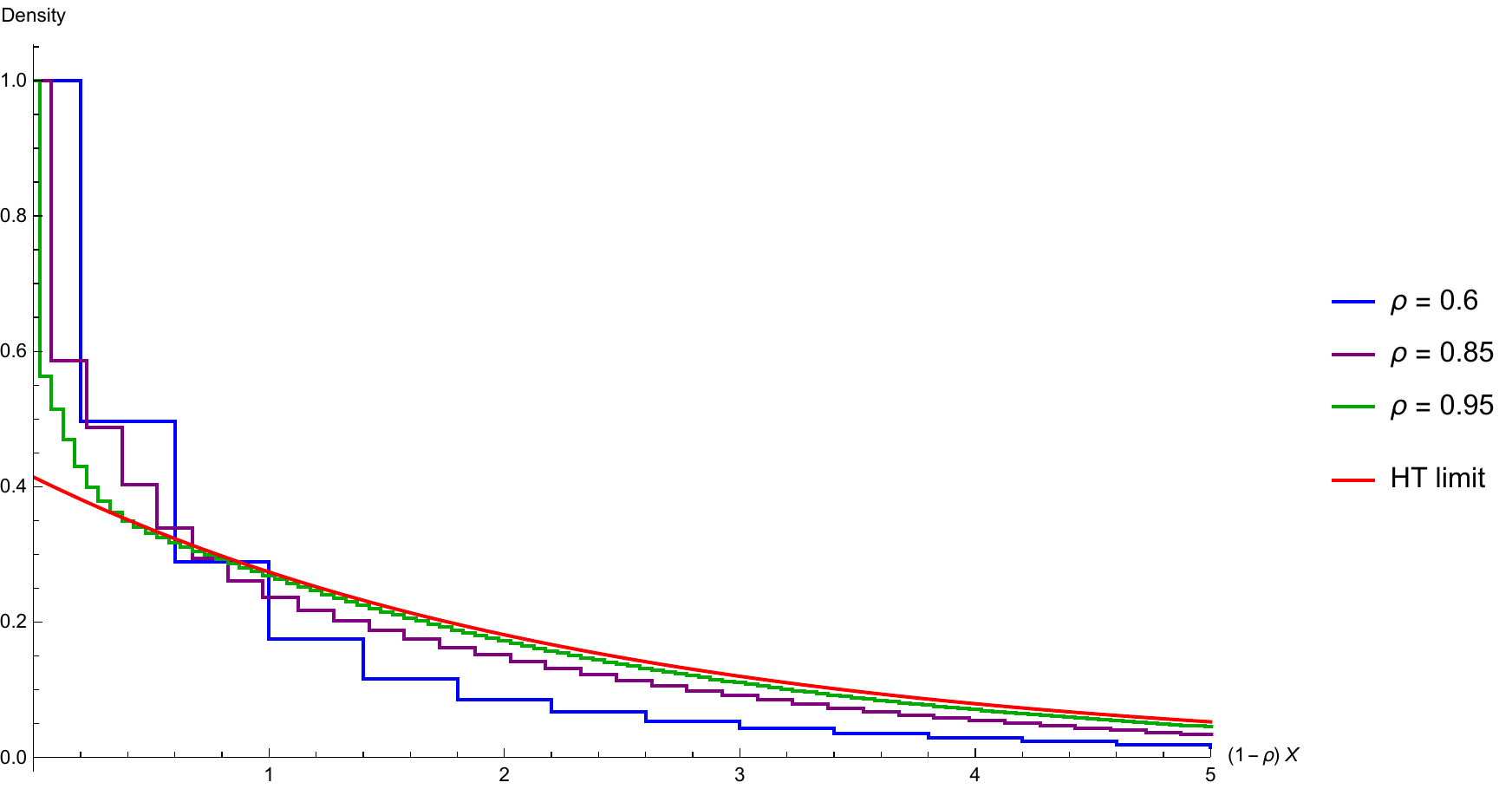}
\end{center}
\caption{The density of the scaled queue length.}
\label{fig:DistXscaledexample}
\end{figure}

Indeed, it can be observed  in Figure \ref{fig:EXscaledexample} that the HT limits of the mean queue lengths in both models, with and without correlated service times, are the same. Note that the \emph{light-traffic} limits, when $\rho\downarrow0$, are also the same. This, however, is caused by the fact that we chose an example with \emph{single} arrivals. In the batch arrival case, the queue-length distributions would also be different in light traffic, due to the correlation between service times of customers inside one batch. It is interesting to see, however, that when  $\rho$ tends to $1$,  the dependence between  subsequent service times no longer influences the mean scaled queue length, and thus the system can be analyzed as an $M/G/1$ queueing system, in this particular example. Furthermore, in Figure \ref{fig:DistXscaledexample}, it can be seen that the density of the scaled queue length converges to the limiting density of an exponential distribution when the traffic intensity $\rho$ approaches $1$.

\section*{Acknowledgments} The research of Abhishek and Rudesindo~N\'u\~nez-Queija  is partly funded by NWO Gravitation project {\sc Networks}, grant number 024.002.003. The authors thank  Onno Boxma (Eindhoven University of Technology) and Michel Mandjes (University of Amsterdam) for helpful discussions.

\bibliographystyle{abbrv}

\end{document}